\documentclass[11pt]{article}
\usepackage{amssymb,amsfonts,amsmath,latexsym,epsf,tikz,url}

\newtheorem{theorem}{Theorem}[section]

\newtheorem{corollary}[theorem]{Corollary}

\newcommand{\proof}{\noindent{\bf Proof.\ }}
\newcommand{\qed}{\hfill $\square$\medskip}

\begin{document}

\title{Distinguishing number and distinguishing index of  	graphs from primary subgraphs} 

\author{ Samaneh Soltani \and
Saeid Alikhani  $^{}$\footnote{Corresponding author}
}

\date{\today}

\maketitle

\begin{center}
Department of Mathematics, Yazd University, 89195-741, Yazd, Iran\\
{\tt s.soltani1979@gmail.com,  alikhani@yazd.ac.ir}
\end{center}

\begin{abstract}
The distinguishing number (index) $D(G)$ ($D'(G)$) of a graph $G$ is the least integer $d$ such that $G$ has an vertex labeling (edge labeling)  with $d$ labels  that is preserved only by a trivial automorphism. Let $G$ be a connected graph constructed from pairwise disjoint connected graphs $G_1,\ldots ,G_k$ by selecting a vertex of $G_1$, a vertex of $G_2$, and identify these two
vertices. Then continue in this manner inductively. We say that $G$ is obtained by point-attaching from $G_1, \ldots ,G_k$ and that $G_i$'s are the primary subgraphs of $G$. 
In this paper, we consider some  particular cases  of these graphs that  are  of importance in chemistry  and study their distinguishing number and index. 
\end{abstract}


\medskip

\section{Introduction}

First, we introduce some notations and terminology which is needed for the the paper. A molecular graph is a simple graph such that its vertices correspond to the atoms and the edges to the bonds of a molecule. Let $G = (V ,E)$ be a graph.  We use the standard graph notation (\cite{Sandi}). In particular, $Aut(G)$ denotes the automorphism group of $G$.  The set of vertices adjacent in $G$ to a vertex of a vertex subset $W\subseteq  V$ is the open neighborhood $N_G(W )$ of $W$. The closed neighborhood  $G[W ]$ also includes all vertices of $W$ itself. In case of a singleton set $W =\{v\}$ we write $N_G(v)$ and $N_G[v]$ instead of $N_G(\{v\})$ and $N_G[\{v\}]$, respectively.  We omit the subscript when the graph $G$ is clear from the context.  The complement of  $N[v]$ in $V(G)$ is  denoted by $\overline{N[v]}$.
 Let $f$ be a mapping from the set $A$ to the set $B$ and let $A' \subseteq A$ and $B' \subseteq B$. If the restriction of $f$ to the set $A'$ is the set $B'$, then we write $f\vert_{A'}=B'$ or $f(A')=B'$.

In theoretical chemistry, molecular structure descriptor, also called topological indices, are used to understand the properties of chemical compounds. The Wiener index is one of the oldest descriptors concerned with the molecular graph \cite{wiener}.   By now there are many different types of such indices for a general graph $G=(V,E)$ (see for example \cite{Deutsch}). Here, apart from the topological index, we are interested in computing the distinguishing number and the distinguishing  index. A labeling of $G$, $\phi : V \rightarrow \{1, 2, \ldots , r\}$, is  $r$-distinguishing,  if no non-trivial  automorphism of $G$ preserves all of the vertex labels. Formally, $\phi$ is $r$-distinguishing if for every non-trivial $\sigma \in Aut(G)$, there exists $x$ in $V$ such that $\phi(x) \neq \phi(x\sigma)$.  The distinguishing number of a graph $G$ is the minimum number $r$ such that $G$ has a labeling that is $r$-distinguishing.   This number has defined by Albertson and Collins \cite{Albert}. Similar to this definition, Kalinowski and Pil\'sniak \cite{Kali1} have defined the distinguishing index $D'(G)$ of $G$ which is  the least integer $d$ such that $G$ has an edge colouring   with $d$ colours that is preserved only by a trivial automorphism.

In this paper, we consider the distinguishing number and the distinguishing index on graphs that contain cut-vertices. Such graphs can be decomposed into subgraphs that we call primary subgraphs. Blocks of graphs are particular examples of primary subgraphs, but a primary subgraph may consist of several blocks. For convenience, the exact definition of these kind of graphs will be given in the next  section.  In Section 2,  the distinguishing number and the distinguishing index of some graphs are computed  from primary subgraphs. In Section 3, we apply the  results of Section 2, in order to obtain the distinguishing number and the distinguishing index  of
families of graphs that are of importance in chemistry.

\section{The distinguishing number (index) of some graphs from primary subgraphs}

\begin{figure}[ht]
	\hspace{3.5cm}
	\includegraphics[width=7cm]{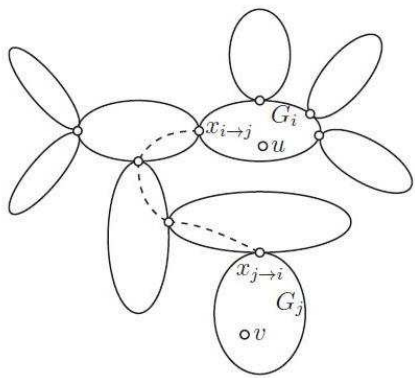}
	\caption{\label{Figure1} Graph $G$ obtained by point-attaching from $G_1,\ldots , G_k$.}
	\end{figure}

Let $G$ be a connected graph constructed from pairwise disjoint connected graphs
$G_1,\ldots ,G_k$ as follows. Select a vertex of $G_1$, a vertex of $G_2$, and identify these two vertices. Then continue in this manner inductively.  Note that the graph $G$ constructed in this way has a tree-like structure, the $G_i$'s being its building stones (see Figure \ref{Figure1}).  Usually  say that $G$ is obtained by point-attaching from $G_1,\ldots , G_k$ and that $G_i$'s are the primary subgraphs of $G$. A particular case of this construction is the decomposition of a connected graph into blocks (see \cite{Deutsch}). 

In this section, we consider some  particular cases  of these graphs  and study their distinguishing number and index.

\begin{figure}[ht]
	\begin{center}
		\includegraphics[width=0.3\textwidth]{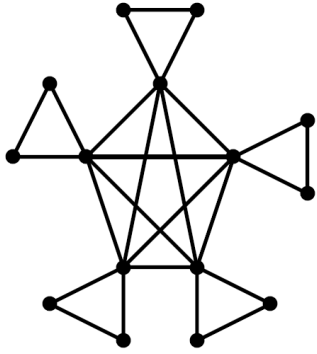}
		\caption{\label{fig1} The graph  $Q(5, 3)$.}
	\end{center} 
\end{figure}

As an example of point-attaching graph,   consider the graph $K_m$ and $m$ copies of  $K_n$. By definition, the graph $Q(m, n)$ is obtained by identifying each vertex of $K_m$ with a vertex of a unique $K_n$. The graph $Q(5, 3)$ is shown in Figure \ref{fig1}.

The following theorems give the distinguishing number and the distinguishing index  of $Q(m,n)$.  

\begin{theorem}
The distinguishing number of $Q(m,n)$ is $$D(Q(m,n))=min\big\{r: r{r\choose n-1}\geqslant m\big\}.$$
\end{theorem}
\proof
We denote the vertices of $K_m$ by $v_1,\ldots , v_m$ and the vertices of corresponding $K_n$ to the vertex $v_i$ by $w_1^{(i)},\ldots, w_{n-1}^{(i)}$ where $i=1,\ldots , m$. In an $r$-distinguishing labeling, each of vertices $w_1^{(i)},\ldots, w_{n-1}^{(i)}$ must have a different labels. Also, each of $n$-ary consisting of a vertex of $K_m$ and $n-1$ vertices of its corresponding $K_n$ must have a different ordered $n$-ary of labels. There are $r{r\choose n-1}$ possible ordered $n$-ary of labels using $r$ labels, hence $D(Q(m,n))=min\{r: r{r\choose n-1}\geqslant m\}$.\qed

\begin{theorem}
The distinguishing index of $Q(m,n)$ is $2$.
\end{theorem}
\proof
We prove the theorem  in three following cases:

Case 1) If $m\geqslant 6$ and $n\geqslant 6$, then we label the edges of $K_m$ and copies of $K_n$ in a distinguishing way with two labels. This labeling is distinguishing, because if $f$ is an automorphism of $Q(m,n)$ preserving  labeling such that it moves the vertices of $K_m$, then with respect to the distinguishing labeling of edges of $K_m$, $f$ dose not preserve the labeling, which is contradiction. So $f$ is the identity automorphism on vertices of $K_m$. With similar argument we can conclude that $f$ is the identity automorphism on vertices of $K_n$, and so on $Q(m,n)$.

\medskip
Case 2) If $m\geqslant 6$ and $n < 6$, then we label the edges of $K_m$ in a distinguishing way with two labels. Since $n < 6$, we can label the edges of every copy of $K_n$ with two labels such that the sets consisting the incident edges to $w^{(i)}_j$, $j=1,\ldots n-1$ have different number of label $2$ for all $i=1,\ldots ,m$. Hence we have a distinguishing labeling for $Q(m,n)$ as prior case.

\medskip
Case 3) If $m < 6$, then we can label the edges of $K_m$ with two labels such that there exist two vertices of $K_m$ have the same number of label $2$ and $1$ in label of their incident edges. Let these two vertices be $v_1$ and $v_2$. We label the edges of $K_n$ corresponding to $v_t$, $t=1,2$ with two labels such that 
\begin{enumerate}
\item[(i)] The sets consisting the incident edges to $w^{(t)}_j$, $j=1,\ldots , n-1$ have different number of label $2$, where $t=1,2$.
\item[(ii)] The number of label $2$ that have been used for the labeling of edges of $K_n$ corresponding to $v_1$ and $v_2$ are distinct.
\end{enumerate}
So this labeling is distinguishing. Using these three cases we have the result.\qed

Now,  we present several constructions of graphs and study their distinguishing number and distinguishing index. These constructions will in turn be used in the next section where chemical applications will be given.  Most of the following constructions have stated in \cite{Deutsch}.

\subsection{Bouquet of graphs}
\begin{figure}[ht]
	\begin{center}
		\includegraphics[width=0.3\textwidth]{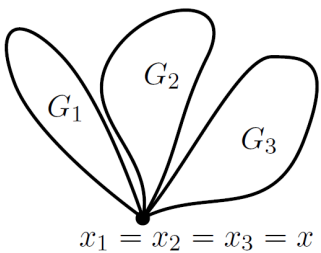}
		\caption{\label{fig2} A bouquet of three graphs.}
	\end{center} 
\end{figure}

Let $G_1,G_2, \ldots ,G_k$ be a finite sequence of pairwise disjoint connected graphs and let
$x_i \in V (G_i)$. By definition, the bouquet $G$ of the graphs $\{G_i\}_{i=1}^k$ with respect to the vertices $\{x_i\}_{i=1}^k$ is obtained by identifying the vertices $x_1, x_2, \ldots ,x_k$ (see Figure \ref{fig2} for $k = 3$).

\begin{theorem}\label{bouquet}
 Let $G$ be the bouquet of the graphs $\{G_i\}_{i=1}^{k}$ with respect to the vertices $\{x_i\}_{i=1}^{k}$, then 
 \begin{enumerate}
 \item[(i)]  $D(G) \leqslant \sum_{i=1}^k D(G_i)$,
 \item[(ii)]  $D'(G) \leqslant \sum_{i=1}^k D'(G_i)$.
 \end{enumerate}
\end{theorem}
\proof 
(i) We label the vertices of graph $G_1$ with labels $\{1,\ldots , D(G_1)\}$ in a distinguishing way. Next we label the vertices of graph $G_j$ ($2 \leqslant j \leqslant k$) except the vertex $x$  with labels $\{(\sum_{i=1}^{j-1} D(G_i) )+1,\ldots , (\sum_{i=1}^{j-1} D(G_i) )+D(G_j)\}$ in a distinguishing way. This labeling is distinguishing, because if $f$ is an automorphism of $G$ preserving the labeling then by the method of labeling we have $f(V(G_i)) = V(G_i)$ where $i=1,\ldots , k$. Since every $G_i$ is labeled distinguishingly, so $f\vert_{V(G_i)}$ is the identity, and so $f$ is the identity automorphism on $G$. We used  $\sum_{i=1}^k D(G_i)$ labels, and hence the result follows.  

(ii) A similar argument yields that $D'(G) \leqslant \sum_{i=1}^k D'(G_i)$. \qed

The bounds of Theorem \ref{bouquet} are sharp. If the graphs $\{G_i\}_{i=1}^{k}$ are the star graphs $\{K_{1,n_i}\}_{i=1}^k$ ($n_i\geqslant 3$) and $\{x_i\}_{i=1}^{k}$ are the central points of $\{K_{1,n_i}\}_{i=1}^k$, respectively, then  the bouquet of $\{K_{1,n_i}\}_{i=1}^k$ with respect to their central points  is the star graph $K_{1,n_1+\cdots + n_k}$. Since the distinguish number and index of the star graph $K_{1,n}$ is $n$, so both bounds of Theorem \ref{bouquet} are sharp.

\subsubsection{Dutch-windmill graphs}
\begin{figure}[ht]
\begin{center}
\includegraphics[width=0.3\textwidth]{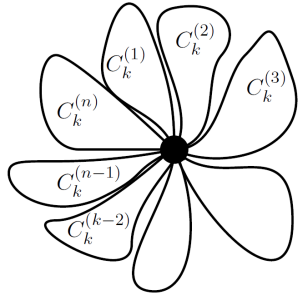}
	\caption{\label{fig01} {\small Dutch windmill Graph $D_n^k$.}}
\end{center}
\end{figure}

Here we consider another kind of point-attaching graphs and study their distinguishing number and distinguishing index.  
The dutch windmill graph $D_n^k$  is the graph obtained by taking $n$, ($n\geqslant 2$) copies of the cycle graph $C_k$, ($k\geqslant 3$) with a vertex in common (see Figure \ref{fig01}). For $k=3$, the graph $D_n^3$ is called the friendship graph and is denoted by $F_n$.  The distinguishing number and the distinguishing index of the friendship graph have been studied in \cite{alikhani}. The following theorem gives the distinguishing number and the distinguishing index of $F_n$. 

\begin{theorem}{\rm \cite{alikhani}}\label{D(friend)}
	\begin{enumerate}
		\item [(i)]
		The distinguishing number of the friendship graph $F_n$  $(n\geq 2)$ is  $$D(F_n)= \lceil \dfrac{1+\sqrt{8n+1}}{2}\rceil .$$  
		\item[(ii)] 
		Let $a_n=1+27n+3\sqrt{81n^2+6n}$. 
		For every $n\geq 2$, $$D'(F_n)=\lceil\frac{1}{3} (a_n)^{\frac{1}{3}}+\frac{1}{3(a_n)^{\frac{1}{3}}}+\frac{1}{3}\rceil.$$
		 \end{enumerate}
\end{theorem}

To obtain the distinguishing number and the distinguishing index of dutch windmill graph $D_n^k$, first  we state and prove the following theorem: 

\begin{theorem}
The order of automorphism group of $ D_n^k $ is $\vert Aut (D_n^k)\vert =n!2^n$.
\end{theorem}
\proof
To obtain  the automorphism group of $D_n^k$, let to denote the central vertex by $w$ and the vertices of $i$th cycle $C_k$ (which we call it a blade) of $D_n^k$ by $V_i$ ($1\leq i\leq n$).  The vertex $w$ should be mapped to itself under  automorphisms of $D_n^k$. In fact every element of the automorphism group of $D_n^k$ is of the form 
\begin{equation*}
h_{\sigma}(v)=\left\{ \begin{array}{ll}
f_1(v)& \textsl{if}~ v\in  V_1\\
\vdots \\
f_n(v)& \textsl{if}~ v\in  V_n
\end{array}\right.
\end{equation*}
where $\sigma \in S_n$. If we denote the vertices of $i$th blade except the central vertex, by $ v^{(i)}_1,\ldots , v^{(i)}_{k-1}$,  then every function $f_i: V_i \rightarrow \sigma (V_i)$ has  one of the following two forms:
\begin{equation*}
\left\{
\begin{array}{l}
v^{(i)}_1\mapsto v^{\sigma (i)}_1\\
v^{(i)}_2\mapsto v^{\sigma (i)}_2\\
\vdots \\
v^{(i)}_{k-1} \mapsto v^{\sigma (i)}_{k-1}
\end{array}\right. \qquad
\left\{
\begin{array}{l}
v^{(i)}_1\mapsto v^{\sigma (i)}_{k-1}\\
v^{(i)}_2\mapsto v^{\sigma (i)}_{k-2}\\
\vdots \\
v^{(i)}_{k-1} \mapsto v^{\sigma (i)}_1
\end{array}\right.
\end{equation*}
Therefore $\vert Aut (D_n^k)\vert =n!2^n$.\qed

\begin{figure}[ht]
\begin{center}
\includegraphics[width=0.7\textwidth]{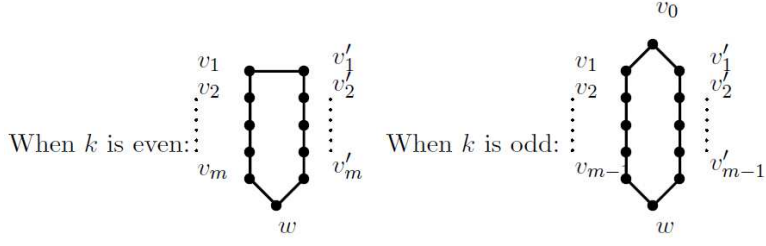}
	\caption{\label{fig0} {\small The considered polygon (or a cycle of size $k$) in the proof of Theorem \ref{thmbasic}}.}
\end{center}
\end{figure}

\begin{theorem}\label{thmbasic}
Let $D_n^k$ be 
dutch windmill graph such that $n\geqslant 2$ and $k\geqslant 3$. Then we have
$D(D_n^k)=min\{r:~\dfrac{r^{k-1}-r^{\lceil \dfrac{k-1}{2} \rceil}}{2}\geqslant n\}$.
\end{theorem}
\proof
If $k$ is odd, then there is a natural number $m$ such that $k=2m+1$. We can consider a blade of $D_n^k$ as Figure \ref{fig0}.
 Let $(x_1^{(i)},x_1'^{(i)},\ldots ,x_m^{(i)},x_m'^{(i)})$ be the label of vertices $(v_1,v_1',\ldots ,v_m,v_m')$ of the $i$th blade where $1\leqslant i \leqslant n$. Suppose that $L=\{(x_1^{(i)},x_1'^{(i)},\ldots ,x_m^{(i)},x_m'^{(i)})\vert ~ 1\leqslant i \leqslant n , x_j^{(i)},x_j'^{(i)}\in \mathbb{N}, 1\leqslant j \leqslant m \}$ is a labeling of the vertices of $D_n^k$ except its central vertex. In an $r$-distinguishing labeling we must have:
\begin{itemize}
\item[(i)] There exists $j\in \{1,\ldots ,m\}$ such that $x_j^{(i)}\neq x_j'^{(i)}$ for all $i\in \{1,\ldots , n\}$.
\item[(ii)] For $i_1\neq i_2$ we must have $(x_1^{(i_1)},x_1'^{(i_1)},\ldots ,x_m^{(i_1)},x_m'^{(i_1)})\neq (x_1^{(i_2)},x_1'^{(i_2)},\ldots ,x_m^{(i_2)},x_m'^{(i_2)})$  and  $(x_1^{(i_1)},x_1'^{(i_1)},\ldots ,x_m^{(i_1)},x_m'^{(i_1)})\neq (x_1'^{(i_2)},x_1^{(i_2)},\ldots ,x_m'^{(i_2)},x_m^{(i_2)})$. 
\end{itemize}
There are $\dfrac{r^{2m}-r^m}{2}$ possible $(2m)$-arrays of labels using $r$ labels satisfying (i) and (ii), hence $D(D_n^k)=min\{r:~\dfrac{r^{2m}-r^m}{2}\geqslant n\}$.

If $k$ is even, then there is a natural number $m$ such that $k=2m$. We can consider a blade of $D_n^k$ as Figure \ref{fig0}. Let $(x_0^{(i)}x_1^{(i)},x_1'^{(i)},\ldots ,x_{m-1}^{(i)},x_{m-1}'^{(i)})$ be the label of vertices $(v_0,v_1,v_1',\ldots ,v_{m-1},v_{m-1}')$ of  $i$th blade where $1\leqslant i \leqslant n$. Suppose that $L=\{(x_0^{(i)},x_1^{(i)},x_1'^{(i)},\ldots ,x_{m-1}^{(i)},x_{m-1}'^{(i)})\vert ~ 1\leqslant i \leqslant n , x_0^{(i)},x_j^{(i)},x_j'^{(i)}\in \mathbb{N}, 1\leqslant j \leqslant m-1 \}$ is a labeling of the vertices of $D_n^k$ except its central vertex. In an $r$-distinguishing labeling we must have:
\begin{itemize}
\item[(i)] There exists $j\in \{1,\ldots ,m-1\}$ such that $x_j^{(i)}\neq x_j'^{(i)}$ for all $i\in \{1,\ldots , n\}$.
\item[(ii)] For $i_1\neq i_2$ we must have 
\begin{align*}
&(x_0^{(i_1)},x_1^{(i_1)},x_1'^{(i_1)},\ldots ,x_{m-1}^{(i_1)},x_{m-1}'^{(i_1)})\neq (x_0^{(i_2)},x_1^{(i_2)},x_1'^{(i_2)},\ldots ,x_{m-1}^{(i_2)},x_{m-1}'^{(i_2)}),\\ &(x_0^{(i_1)},x_1^{(i_1)},x_1'^{(i_1)},\ldots ,x_{m-1}^{(i_1)},x_{m-1}'^{(i_1)})\neq (x_0^{(i_2)},x_1'^{(i_2)},x_1^{(i_2)},\ldots ,x_{m-1}'^{(i_2)},x_{m-1}^{(i_2)}).
\end{align*}
\end{itemize}
There are $\dfrac{r^{2m-1}-r^m}{2}$ possible $(2m-1)$-arrays of labels using $r$ labels satisfying (i) and (ii) ($r$ choices for $x_0$ and $\dfrac{r^{2(m-1)}-r^{m-1}}{2}$ choices for $x_1^{(i_1)},x_1'^{(i_1)},\ldots ,x_{m-1}^{(i_1)},x_{m-1}'^{(i_1)}$), hence $D(D_n^k)=min\{r:~\dfrac{r^{2m-1}-r^m}{2}\geqslant n\}$.\qed

\begin{corollary}\label{nice}
Let $D_n^k$ be the dutch windmill graph such that $n\geqslant 2$ and $k\geqslant 3$. If $k=2m+1$, then $D(D_n^k)=\lceil \sqrt[m]{\dfrac{1+\sqrt{8n+1}}{2}}\rceil$.
\end{corollary}
\proof
It is easy to see that $min\{r:~\dfrac{r^{2m}-r^m}{2}\geqslant n\}=\lceil \sqrt[m]{\dfrac{1+\sqrt{8n+1}}{2}}\rceil$. So the result follows from Theorem \ref{thmbasic}.
\qed

The following theorem implies that to study the distinguishing index of $D_n^k$, it suffices to study its distinguishing number and vice versa. 
\begin{theorem}\label{thmindex}
Let $D_n^k$ be the dutch windmill graph such that $n\geqslant 2$ and $k\geqslant 3$. Then 
$D'(D_n^k)=D(D_n^{k+1})$.
\end{theorem}
\proof
Since the effect of every automorphism of $D_n^{k+1}$ on its non-central vertices is exactly the same as the effect of an automorphism of $D_n^k$ on its edges and vice versa, so if we consider the non-central vertices of $D_n^{k+1}$ as the edges of $D_n^k$, then we have $D'(D_n^k)=D(D_n^{k+1})$.
\qed

\subsection{Circuit of graphs}
\begin{figure}[ht]
	\begin{center}
		\includegraphics[width=0.34\textwidth]{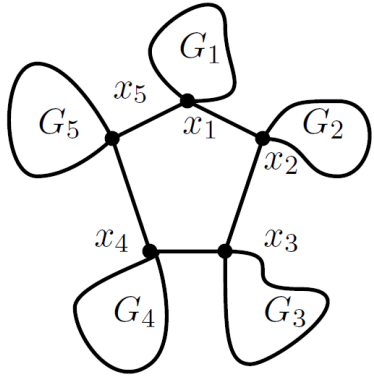}
		\caption{\label{fig3} A circuit of five graphs.}
	\end{center} 
\end{figure}

Let $G_1,G_2, \ldots ,G_k$ be a finite sequence of pairwise disjoint connected graphs and let
$x_i \in  V (G_i)$. By definition, the circuit $G$ of the graphs $\{G_i\}_{i=1}^k$ with respect to the vertices $\{x_i\}_{i=1}^k$ is obtained by identifying the vertex $x_i$ of the graph $G_i$ with the $i$-th vertex of the cycle graph $C_k$ (\cite{Deutsch}). See Figure \ref{fig3} for $k = 5$.

\begin{theorem}
Let $G$ be circuit graph of the graphs $\{G_i\}_{i=1}^k$ with respect to the vertices $\{x_i\}_{i=1}^k$. Then
\begin{enumerate}
\item[(i)] $D(G)\leqslant max\big\{ max \{D(G_i)\}_{i=1}^k ,D(C_k)\big\}$,
\item[(ii)] $D'(G)\leqslant max\big\{ max \{D'(G_i)\}_{i=1}^k ,D'(C_k)\big\}$.
\end{enumerate}
\end{theorem}
\proof
(i) We label the vertices of $C_k$ with the labels $\{1,\ldots , D(C_k)\}$ and vertices of every $G_i$ ($1\leqslant i \leqslant k$) with the labels $\{1,\ldots , D(G_i)\}$ in a distinguishing way, respectively. This labeling is distinguishing for $G$, because if $f$ is an automorphism of $G$ preserving the labeling, then we have two following cases:
\begin{enumerate}
\item[(a)] If $f\vert_{V(C_k)}=V(C_k)$, then for all $i$ we have $f\vert_{V(G_i)}=V(G_i)$.  Since we labeled $C_k$ distinguishingly, $f\vert_{V(C_k)}$ is the identity automorphism. In this case $f$ is the identity automorphism on $G$, because each of $G_i$ is labeled in a distinguishing way. We used $max\{ max \{D(G_i)\}_{i=1}^k ,D(C_k)\}$ labels, and so the result follows.
 \item[(b)] Suppose that there exists the vertex $x$ of $C_k$ such that for some $i$, $f(x)=y$ where $y \in V(G_i)\setminus V(C_k)$ and $x\notin V(G_i)$, then $f(V(C_k)) \subseteq V(G_i)$, i.e., $G_i$ is contains a copy of $C_k$. The label of vertex $y$ can be $l\in \{1,\ldots , D(G_i)\}$ and label $x$ can be $l'\in \{1,2,3\}$. By assigning  two different labels to $x$ and $y$ we get a distinguishing labeling
\end{enumerate}

(ii) A similar argument yields that $D'(G)\leqslant max\big\{ max \{D'(G_i)\}_{i=1}^k ,D'(C_k)\big\}$.\qed

\subsection{Chain of graphs}
\begin{figure}[ht]
	\begin{center}
		\includegraphics[width=0.7\textwidth]{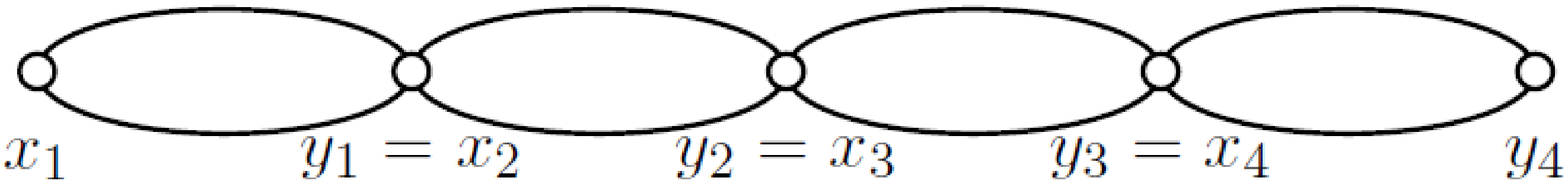}
		\caption{\label{fig4} A chain of graphs.}
	\end{center} 
\end{figure}

Let $G_1,G_2,\ldots ,G_k$ be a finite sequence of pairwise disjoint connected graphs and let
$x_i, y_i \in V (G_i)$. By definition (see \cite{Mansour,Schork}) the chain $G$ of the graphs $\{G_i\}_{i=1}^k$ with respect to the vertices $\{x_i, y_i\}_{i=1}^k$ is obtained by identifying the vertex $y_i$ with the vertex $x_{i+1}$ for $i \in [k - 1]$ (also see \cite{Deutsch}). See Figure \ref{fig4} for $k = 4$.
\begin{theorem}
Let $G$ be the chain  of the graphs $\{G_i\}_{i=1}^k$ with respect to the vertices $\{x_i, y_i\}_{i=1}^k$. Then
\begin{enumerate}
\item[(i)] $D(G) \leqslant max \big\{max\{D(G_i)\}_{i=1}^k, max \{deg_G x_i \big\}_{i=2}^k\}$,
\item[(ii)] $D'(G) \leqslant max \big\{max\{D'(G_i)\}_{i=1}^k, max \{deg_G x_i \big\}_{i=2}^k\}$.
\end{enumerate} 
\end{theorem}
\proof
 The shortest path between $x_1$ and $y_k$ is made by connecting the shortest paths between $x_i$ and $y_i$ for $i=1,\ldots , k$. If $f$ is an automorphism of $G$, then we have the two following cases:
 \begin{enumerate}
\item[$(a)$]  There exists $i$ ($1\leqslant i \leqslant k$) such that $f(x_i)\neq x_i$. Thus the shortest path between $x_1$ and $y_k$ is not fixed under $f$, and so $f(x_1)=y_k$ or $f(y_k)=x_1$.
\item[$(b)$]  For all $i$, $f(x_i) = x_i$, and so $f(N_G(x_i))= N_G(x_i)$.
\end{enumerate} 
\begin{enumerate}
\item[(i)] Now we want to present a distinguishing vertex labeling for $G$. For this purpose we label the vertex $x_1$ with label $1$ and the vertex $y_k$ with label $2$, next we label all vertices adjacent to $x_i$, with $deg_G x_i$ different labels where $i = 2,\ldots , k$. The remaining  vertices of every $G_i$ is labeled with labels  $\{1, \ldots , D(G_i)\}$ in a distinguishing way, respectively. This labeling is distinguishing, because if $f$ is an automorphism of $G$ preserving the labeling then the two following cases may occur:
 \begin{enumerate}
\item[$(a')$]  There exists $i$ ($1\leqslant i \leqslant k$) such that $f(x_i)\neq x_i$. Thus by (a)  we have $f(x_1)=y_k$ or $f(y_k)=x_1$. Since we label the vertices  $x_1$ and $y_k$ with two different labels, so this case can not occur.
\item[$(b')$]  For all $i$, $f(x_i) = x_i$, and so by (b) we have $f(N_G(x_i))= N_G(x_i)$. Since adjacent vertices to every $x_i$ are labeled differently, so $f\vert_{N_G(x_i)}$ is the identity automorphism. On the other hand since $f(N_G[x_i])= N_G[x_i]$, hence $f\vert_{V(G_i)}=V(G_i)$, and thus $f\vert_{V(G_i)}$ is the identity automorphism, because we labeled the vertices in $\overline{N_{G_i}(x_i)}$ distinguishingly. Therefore $f$ is the identity automorphism on $G$.
\end{enumerate}
Since we used $max \{max\{D(G_i)\}_{i=1}^k, max \{deg_G x_i \}_{i=2}^k\}$ labels, the Part (i) follows.

\item[(ii)] First we label all edges incident to $x_1$ with label $1$, and all edges incident to $y_k$ with label $2$.  Next we label all edges incident to $x_i$, with $deg_G x_i$ different labels where $i = 2,\ldots , k$. The remaining  edges of every $G_i$ is labeled with labels  $\{1, \ldots , D'(G_i)\}$ in a distinguishing way, respectively. As Case (i) we can prove this labeling is distinguishing, and that $D'(G) \leqslant max \{max\{D'(G_i)\}_{i=1}^k, max \{deg_G x_i \}_{i=2}^k\}$.\qed
\end{enumerate}

\subsection{Link of graphs}
\begin{figure}[ht]
	\begin{center}
		\includegraphics[width=0.7\textwidth]{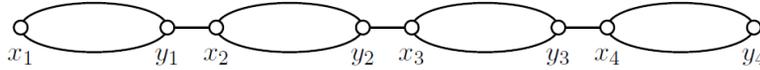}
		\caption{\label{fig6} A link of graphs.}
	\end{center} 
\end{figure}

Let $G_1,G_2,\ldots ,G_k$ be a finite sequence of pairwise disjoint connected graphs and let
$x_i, y_i \in V (G_i)$. By definition (see \cite{Ghorbani}), the link $G$ of the graphs $\{G_i\}_{i=1}^k$ with respect to the vertices $\{x_i, y_i\}_{i=1}^k$  is obtained by joining by an edge the vertex $y_i$ of $G_i$ with the vertex $x_{i+1}$ of $G_{i+1}$ for all $i = 1,2,\ldots, k - 1$ (see Figure \ref{fig6} for $k = 4$).

\begin{theorem}
Let $G$ be the link   of the graphs $\{G_i\}_{i=1}^k$ with respect to the vertices $\{x_i, y_i\}_{i=1}^k$. Then
\begin{enumerate}
\item[(i)] $D(G)\leqslant max\{D(G_i)\}_{i=1}^k$,
\item[(ii)] $D'(G)\leqslant max\{D'(G_i)\}_{i=1}^k$.
\end{enumerate}
\end{theorem}
\proof  The shortest path between $x_1$ and $y_k$ is made by connecting the shortest paths between $x_i$ and $y_i$ for $i=1,\ldots , k$, altogether the edges $y_ix_{i+1}$ where $i=1,\ldots, k-1$. If $f$ is an automorphism of $G$, then we have the two following cases:
 \begin{enumerate}
\item[$(a)$]  There exists $i$ ($1\leqslant i \leqslant k$) such that $f(x_i)\neq x_i$. Thus the shortest path between $x_1$ and $y_k$ is not fixed under $f$, and so $f(x_1)=y_k$ or $f(y_k)=x_1$.
\item[$(b)$]  For all $i$, $f(x_i) = x_i$, and so $f(y_i)= y_i$, and hence $f(V(G_i))=V(G_i)$.
\end{enumerate} 

\begin{enumerate}
\item[(i)] Now we want to present a distinguishing vertex labeling for $G$. For this purpose we label the vertex $x_1$ with label $1$ and the vertex $y_k$ with label $2$. The   vertices of every $G_i$ is labeled with labels  $\{1, \ldots , D(G_i)\}$ in a distinguishing way, respectively. This labeling is distinguishing, because if $f$ is an automorphism of $G$ preserving the labeling then the two following cases may occur:
 \begin{enumerate}
\item[$(a')$]  There exists $i$ ($1\leqslant i \leqslant k$) such that $f(x_i)\neq x_i$. Thus by (a)  we have $f(x_1)=y_k$ or $f(y_k)=x_1$. Since we labeled the vertices  $x_1$ and $y_k$ with two different labels, so this case can not occur.
\item[$(b')$]  For all $i$, $f(x_i) = x_i$, and so by (b) we have $f(y_i)= y_i$ and $f(V(G_i))=V(G_i)$. Hence $f\vert_{V(G_i)}$ is the identity automorphism, because we labeled the vertices of each $G_i$ distinguishingly. Therefore $f$ is the identity automorphism on $G$.
\end{enumerate}
Since we used $max\{D(G_i)\}_{i=1}^k$ labels, the Part (i) follows.

\item[(ii)] First we label all edges incident to $x_1$ with label $1$, and all edges incident to $y_k$ with label $2$.  The edges of every $G_i$ is labeled with labels  $\{1, \ldots , D'(G_i)\}$ in a distinguishing way, respectively. As Case (i) we can prove this labeling is distinguishing, and that $D'(G) \leqslant max\{D'(G_i)\}_{i=1}^k$.\qed
\end{enumerate}

\section{Distinguishing labeling of graphs that are of importance in chemistry }

In this section, we apply the previous results in order to obtain the distinguishing number and the distinguishing index  of
families of graphs that are of importance in chemistry.

\subsubsection{Spiro-chains}
\begin{figure}[ht]
	\begin{center}
		\includegraphics[width=0.45\textwidth]{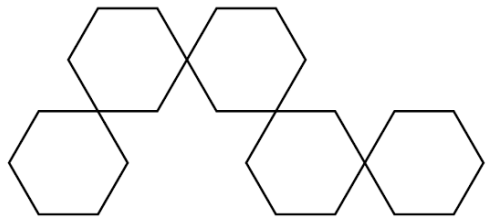}
		\caption{\label{fig5} Spiro-chain $S_{6,2,5}$.}
	\end{center} 
\end{figure}

Spiro-chains are defined in \cite{Diudea} page $114$. Making use of the concept of chain of graphs,
a spiro-chain can be defined as a chain of cycles. We denote by $S_{q,h,k}$ the chain of
$k$ cycles $C_q$ in which the distance between two consecutive contact vertices is $h$ (see
$S_{6,2,5}$ in Figure \ref{fig5}). 

\begin{theorem}\label{DSpiro}
	The distinguishing number of spiro-chain  $S_{q,h,k}$ is $2$, except $D(S_{3,1,2})=3$.
\end{theorem}
\proof
Since spiro-chains are as a chain of cycles, we follow the notation of attached vertices   as denoted in the chain of graphs in Figure \ref{fig4}. We assign the vertex $x_1$, label $1$, and the vertices $y_1,\ldots,y_k$, label $2$. Next we assign non-labeled vertices on path of length $h$ between $x_i$ and $y_i$, label $1$ and assign the non-labeled vertices on path of length $k-h$ between $x_i$ and $y_i$, the label $2$, where $1 \leqslant i \leqslant k$. This labeling is distinguishing, because with respect to the label of $x_1$ and $y_k$, every $C_k$ mapped to itself. Now regarding to the label of vertices of each $C_k$, we can conclude that the identity automorphism is the only automorphism of $G$ preserving the labeling. So $D(S_{q,h,k})=2$. 
It can be seen that  $S_{3,1,2}$ is friendship graph $F_2$, so $D(S_{3,1,2})=3$ by Theorem \ref{D(friend)}.\qed

\begin{theorem}
	The distinguishing index of spiro-chain $S_{q,h,k}$ is $2$.
\end{theorem}
\proof
We label the two edges incident to $x_1$ with label $1$, and the two edges incident to $y_k$ with label $2$. Next we assign non-labeled edges on path of length $h$ between $x_i$ and $y_i$, label $1$ and assign the non-labeled edges on path of length $k-h$ between $x_i$ and $y_i$, label $2$  where $1 \leqslant i \leqslant k$. Similar the proof of Theorem \ref{DSpiro}, this labeling is distinguishing. So $D'(S_{q,h,k})=2$. \qed

\subsubsection{Polyphenylenes}
\begin{figure}[ht]
	\begin{center}
		\includegraphics[width=0.7\textwidth]{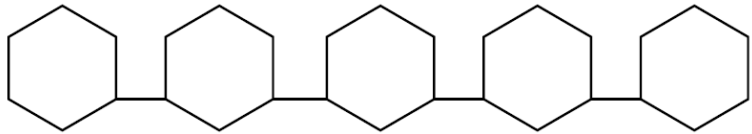}
		\caption{\label{fig7} Polyphenylenes $L_{6,2,5}$.}
	\end{center} 
\end{figure}

Similar  to the  definition of the spiro-chain $S_{q,h,k}$, we can define the graph
$L_{q,h,k}$ as the link of $k$ cycles $C_q$ in which the distance between the two contact
vertices in the same cycle is $h$. (See Figure \ref{fig7} for $L_{6,2,5}$).

 The following theorem shows that polyphenylenes can be distinguished by
 two labels. 
\begin{theorem} 
	\begin{enumerate}
		\item[(i)] The distinguishing number of polyphenylene $L_{q,h,k}$ is $2$.
		
		\item[(ii)] The distinguishing index of polyphenylene $L_{q,h,k}$  is $2$.
		
		\end{enumerate}
	\end{theorem}
\proof
\begin{enumerate}
 \item[(i)]  The proof is exactly similar to the proof of Theorem \ref{DSpiro}.
 
 \item[(ii)] Since polyphenylenes are as a link  of cycles, we follow the notation of attached vertices   as denoted in the link of graphs in Figure \ref{fig6}. We label the  edge $x_1y_1$ with label $1$, and the  edges $x_2y_2,\ldots , x_ky_k$ with label $2$. Next we assign the non-labeled edges of every $C_k$ on path of length $h$ between $x_i$ and $y_i$, label $1$ where $1 \leqslant i \leqslant k$, and we label the rest of edges of $C_k$ with label $2$. It can be seen that this labeling is distinguishing. So $D'(L_{q,h,k})=2$. \qed
\end{enumerate}

\subsection{Nanostar dendrimers}
\begin{figure}[ht]
	\begin{center}
		\includegraphics[width=0.4\textwidth]{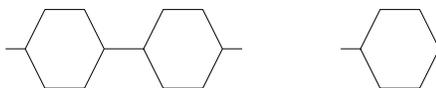}
		\caption{\label{fig8} The graphs $F$ and $G_1$, respectively.}
	\end{center} 
\end{figure}

\begin{figure}[ht]
	\begin{center}
		\includegraphics[width=0.5\textwidth]{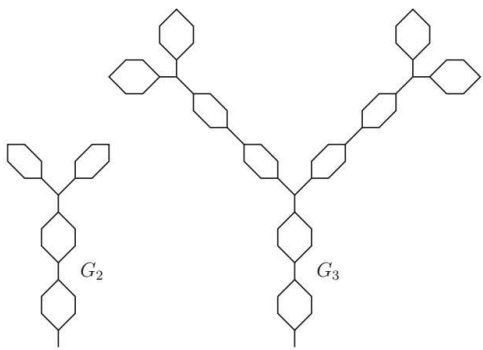}
		\caption{\label{G2} The graphs $G_2$ and $G_3$, respectively.}
	\end{center} 
\end{figure}

\begin{figure}[ht]
	\begin{center}
		\includegraphics[height=11cm, width=9.8cm]{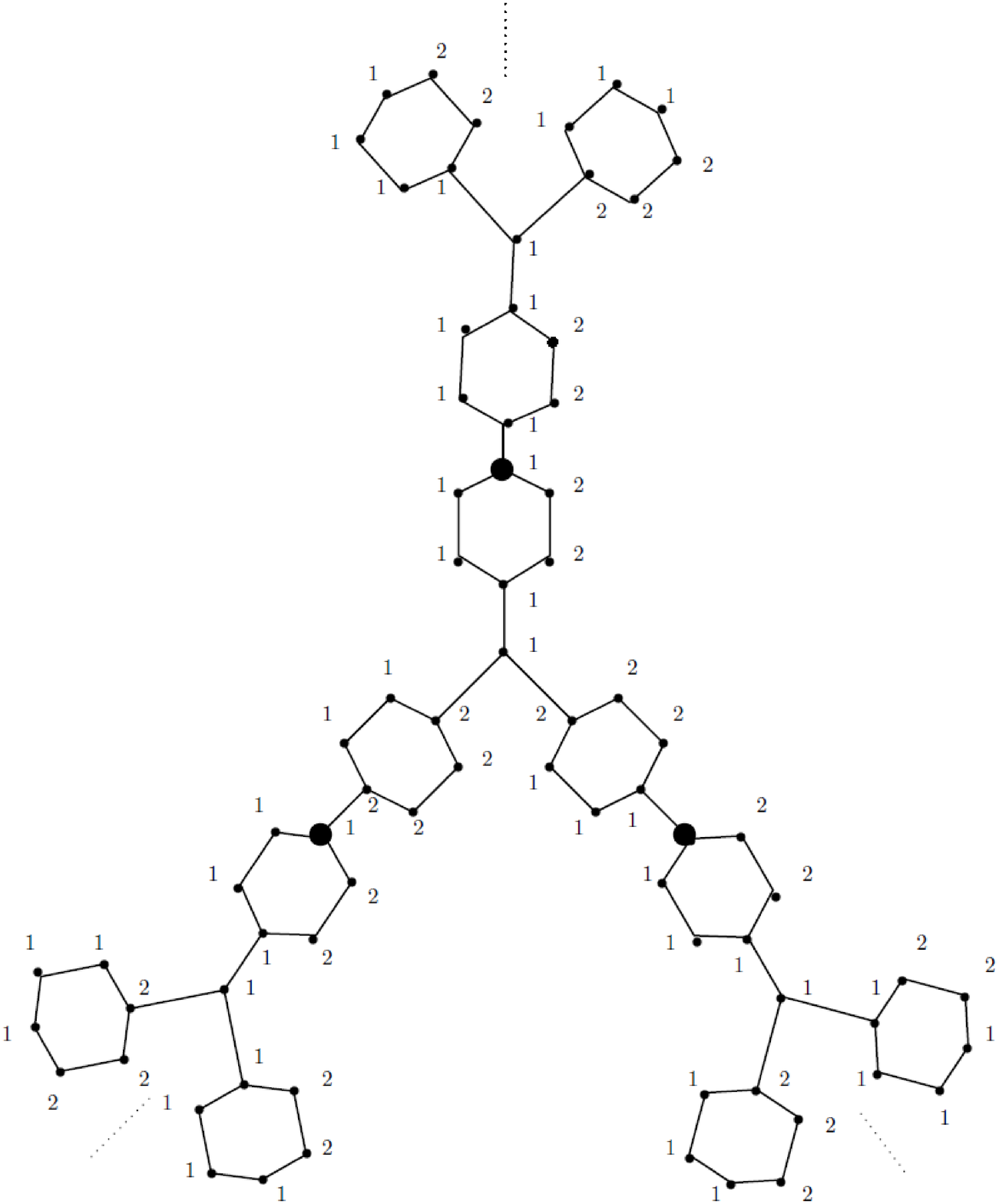}
		\caption{\label{fig9} The $2$-distinguishing labeling of vertices of $ND_n$.}
	\end{center} 
\end{figure}

Dendrimers are large and complex molecules with well taylored
chemical structures. These are key molecules in nanotechnology and can be put to good use. 
We intend to derive the distinguishing number and the distinguishing index of  the nanostar dendrimer $ND_k$ defined  in \cite{17}.  In order to define $ND_k$, we follow \cite{Deutsch}. First we define recursively an auxiliary family of rooted dendrimers $G_k$ $(k \geq 1)$. We need 
a fixed graph $F$ defined in Figure \ref{fig8}; we consider one of its endpoint to be the root of $F$. The graph $G_1$ is defined in Figure \ref{fig8}, the leaf being its root. Now we define $G_k$ $(k\geq 2)$ as the bouquet of the following 3 graphs: $G_{k-1}, G_{k-1}$, and $F$ with respect to their roots; the
root of $G_k$ is taken to be its unique leaf (see $G_2$ and $G_3$ in Figure \ref{G2}). Finally, we define $ND_k$ $(k \geq  1)$ as the bouquet of 3 copies of $G_k$ with respect to their roots.  See a nanostar dendrimer $ND_n$ depicted in Figure \ref{fig9}.

One can see that the nanostar dendrimer graph is a partial cube. In this subsection  we  compute the distinguishing number and the distinguishing index of this infinite class of dendrimers. 

\begin{theorem}
The distinguishing number and index of nanostar dendrimer graph is $2$. 
\end{theorem}
\proof
 Since the nanostar dendrimer graph is symmetric, so $D(ND_n) > 1$. In Figure \ref{fig9} we presented a $2$-labeling of vertices of $ND_n$. Considering the symmetries (automorphisms) of $ND_n$, it can follow  that the labeling is distinguishing.  A similar argument also yields $D'(ND_n)=2$.
  \qed

\end{document}